\documentclass{amsart}
\usepackage{amssymb}

\newtheorem{thm}{Theorem}
\newtheorem{prop}
[thm]{Proposition}
\newtheorem{lem}
[thm]{Lemma}

\def\geq{\geqslant}

\begin{document}

\title{Kleinian groups which are
almost fuchsian}

\author{Francis Bonahon}

\address {F. Bonahon, Department
of Mathematics,  University of
Southern California, Los Angeles,
CA~90089-1113, U.S.A.}
\email{fbonahon@math.usc.edu}
\urladdr{http://math.usc.edu/\~{}%
fbonahon}

\thanks{This work was partially
supported by grant DMS-0103511
from the National Science
Foundation.}

\begin{abstract} We consider the
space of all quasifuchsian metrics
on the product of a surface with
the real line. We show that, in a
neighborhood of the submanifold
consisting of fuchsian metrics,
every non-fuchsian metric is
completely determined by the
bending data of its convex core. 
\end{abstract}

\date{\today}
\maketitle

Let $S$ be a surface of finite
topological type, obtained by
removing finitely many points from
a compact surface without
boundary, and with negative Euler
characteristic. We consider
complete hyperbolic metrics on
the product $S\times
\left]-\infty,
\infty \right[$. 

The simplest ones are
the fuchsian metrics defined as
follows. Because of our
hypothesis that the Euler
characteristic of $S$ is
negative, $S$ admits a finite
area hyperbolic metric, for which
$S$ is isometric to the quotient
of the hyperbolic plane $\mathbb
H^2$ by a discrete group $\Gamma$
of isometries. The group $\Gamma$
uniquely extends to a group of
isometries of the hyperbolic
3--space $\mathbb H^3$ respecting
the transverse orientation of
$\mathbb H^2 \subset \mathbb
H^3$, for which the quotient
$\mathbb H^3/\Gamma$ has a
natural identification with 
$S\times \left]-\infty,
\infty \right[$. A \emph{fuchsian
metric} is any metric on 
$S\times \left]-\infty,
\infty \right[$ obtained
in this way. Note that the image
of
$\mathbb H^2$ in $\mathbb H^3$
provides in this case a totally
geodesic surface in $S\times
\left]-\infty,
\infty \right[$, isometric to the
original metric on $S$.

These examples can be perturbed
to more complex hyperbolic
metrics on $S\times \left]-\infty,
\infty \right[$. See for instance
\cite{Thu1}\cite{Mas}. A
\emph{quasifuchsian metric} on 
$S\times \left]-\infty,
\infty \right[$ is one which is
obtained by quasi-conformal
deformation of a fuchsian metric.
Equivalently, a quasifuchsian
metric is a geometrically finite
hyperbolic metric on
$S\times \left]-\infty,
\infty \right[$ whose cusps
exactly correspond to the ends of
$S$. These also correspond to the
interior points in the space of
all hyperbolic metrics on 
$S\times \left]-\infty,
\infty \right[$ for which the ends
of $S$ are parabolic
\cite{Mar}\cite{Sul}.

If $m$ is a quasifuchsian metric
on $S\times \left]-\infty,
\infty \right[$, the totally geodesic copy of
$S$ which occurred in the fuchsian
case is replaced by the
\emph{convex core} $C(m)$,
defined as the smallest non-empty
closed $m$--convex subset of 
$S\times \left]-\infty,
\infty \right[$. If $m$ is not
fuchsian, $C(m)$ is
3--dimensional and its boundary
consists of two copies of $S$,
each facing an end of $S\times
\left]-\infty,
\infty \right[$. The geometry of
$\partial C(m)$ was investigated
by Thurston
\cite{Thu1}; see
also \cite{EpsMar}.  The
component of $\partial C(m)$ that
faces the end
$S\times\{+\infty\}$ is a
pleated surface, totally geodesic
almost everywhere, but bent along
a family of simple geodesics;
this bending is described and
quantified by a measured geodesic
lamination $\beta^+(m)$ on $S$.
Similarly, the bending of the
negative component of $\partial
C(m)$, namely the one facing
$S\times\{-\infty\}$, is
determined by a measured geodesic
lamination $\beta^-(m)$. 

If $\mathcal Q(S)$ denotes the
space of isotopy classes of
quasifuchsian metrics on $S\times \left]-\infty,
\infty \right[$ and if $\mathcal
{ML}(S)$ is the space of measured
geodesic laminations on $S$, the
rule $m \mapsto \left(
\beta^+(m), \beta^-(m) \right)$
defines a map $\beta: \mathcal
Q(S) \rightarrow \mathcal {ML}(S)
\times \mathcal {ML}(S)$. By
definition, $\beta(m)=(0,0)$ is
the metric $m$ is fuchsian, in
which case the convex core $C(m)$
is just a totally $m$--geodesic
copy of $S$. The image
$\beta(m)\in \mathcal {ML}(S)^2$,
interpreted as a measured
geodesic lamination on two copies
of $S$, is the
\emph{bending measured geodesic
lamination} of the quasifuchsian
metric $m$. 

The space $\mathcal Q(S)$ is a
manifold of dimension $2\theta$
where, if $\chi(S)$ is the Euler
characteristic of $S$ and
$p$ is its number of ends, 
$\theta=-3\chi(S)-p\geq 0$. It
contains the space $\mathcal
F(S)$ as a proper submanifold of
dimension $\theta$. To some
extend, the bending measured
lamination $\beta(m)$ measures
how far the metric $m \in
\mathcal Q(S)$ is from being
fuchsian. Finally, recall that
the space $\mathcal {ML}(S)$ of
measured geodesic laminations is
a piecewise linear manifold of
dimension $\theta$. 

Thurston conjectured that the
restriction of the bending map to
$\beta$ could be used to
parametrize $\mathcal Q(S) -
\mathcal F(S)$, namely that it
induces a homeomorphism between
$\mathcal Q(S) - \mathcal F(S)$
and an open subset of $\mathcal
{ML}(S)^2$. The image of $\beta$
was determined in \cite{BonOta}. 

The goal of the current paper is
to prove Thurston's conjecture on
a neighborhood of the space of
Fuchsian metrics. 

\begin{thm} 
\label{thm:MainThm}
There exists
an open neighborhood $V$ of the
fuchsian submanifold $\mathcal
F(S)$ in $\mathcal Q(S)$
such that the
bending map $\beta: \mathcal Q(S)
\rightarrow \mathcal{ML}(S)^2$
induces a homeomorphism between
$V-\mathcal F(S)$ and its image.  
\end{thm}

There are well-known restrictions
for $(\mu, \nu)\in \mathcal
{ML}(S)^2$ to be in the image of
$\beta$; see for instance
\cite{BonOta}. In particular, if
$(\mu, \nu)\not = (0,0)$ is the
bending measured lamination of
some $m\in \mathcal Q(S)$, then
the measured geodesic laminations
$\mu$ and $\nu$ must \emph{fill
up} the surface $S$, in the sense
that  every non-trivial measured
geodesic lamination has non-zero
geometric intersection number
with at least one of $\mu$,
$\nu$. This is equivalent to the
condition that every component of
$S-\mu \cup
\nu$ is, either a topological disk
 bounded by the union of
finitely many geodesic arcs, or
a topological annulus bounded on
one side by the union of finitely
many geodesic arcs and going to a
cusp on the other side. 

Let $\mathcal {FML}(S)$ denote
the open subset of
$\mathcal{ML}(S)^2$ consisting of
those $(\mu, \nu)$ where $\mu$
and $\nu$ fill up
$S$. Note that $\mathcal
{FML}(S)$ is endowed with an
action of $\mathbb R^+$, defined
by $t(\mu, \nu)=(t\mu, t\nu)$,
which decomposes $\mathcal
{FML}(S)$ as the union of
pairwise disjoint rays (=orbits)
$\left]0,\infty \right[(\mu,
\nu)$.

\begin{thm}
\label{thm:LocalImageBeta}
 In
Theorem~\ref{thm:MainThm}, the
neighborhood $V$ of $\mathcal
F(S)$ and its
image $U=\beta(V)$ can be chosen
so that $U-\{(0,0)\}$ is an open
subset of
$\mathcal {FML}(S)$ which
intersects each ray 
$\left]0,\infty \right[(\mu,
\nu)$ in an interval 
$\left]0,\varepsilon_{\mu\nu}
\right[(\mu,
\nu)$.
\end{thm}

Theorems~\ref{thm:MainThm} and
\ref{thm:LocalImageBeta} are
proved later as
Theorem~\ref{thm:HomeoNoCheck}.
The main idea of the proof is
to construct an inverse
$\beta^{-1}: U-\{(0,0)\}
\rightarrow V-\mathcal F(S)$ ,
and splits into two steps: an
infinitesimal part, and a
transversality argument based on
the infinitesimal part. The
infinitesimal part is now
relatively classical; see for
instance \cite{Ser1}. There are
restrictions on which
 bending data can be
realized by an infinitesimal
deformation of $m_0\in \mathcal
F(S)$. Through the complex
structure of $\mathcal Q(S)$,
where multiplication by $\mathrm
i=\sqrt{-1}$ converts bending to
shearing, these restrictions can
be expressed in a purely
2--dimensional context; see
Section~\ref{sect:NecCond}. The
main part of the proof is to show
by a transversality argument that
any infinitesimal bending data
can actually be realized by a
deformation. The only significant
idea of the paper is to apply the
transversality argument, not in
$\mathcal Q(S)$ where the
necessary hypotheses are not
realized, but in the
manifold-with-boundary $\check
{\mathcal Q}(S)$ obtained by
blowing up $\mathcal Q(S)$ along
the fuchsian submanifold
$\mathcal F(S)$. 

Switching to the blow-up
manifold $\check{\mathcal Q}(S)$
actually provides a better
understanding of the
restriction of $\beta$ on a
neighborhood of $\mathcal F(S)$
and of its inverse. See
Theorem~\ref {them:HomeoInCheck}
for a precise statement.

\medskip
\noindent\textbf
{Acknowledgements:}
There is strong evidence that the
content of this paper has been
known to Bill Thurston for about
twenty years. I personally
learned much of the material of
Section~\ref{sect:NecCond} many
years ago from David Epstein
who, I believe, had learned it
from Thurston. I am also
grateful to Pete Storm and Dick
Canary for helping me to clarify
my thoughts. 

\section{The earthquake section}
\label{sect:EarthSec}

We consider the
\emph{Teichm\"uller space}
$\mathcal T(S)$, namely the space
of isotopy classes of finite area
complete hyperbolic metrics on
the surface $S$. Recall that
$\mathcal T(S)$ is diffeomorphic
to $\mathbb R^\theta$, where 
$\theta=-3\chi(S)-p\geq 0$
if $\chi(S)$ is the Euler
characteristic of $S$ and
$p$ is its number of ends.

A standard
deformation of a metric $m\in
\mathcal T(S)$ is the \emph{left
earthquake} $E^\mu_m \in \mathcal
T(S)$ along the measured geodesic
lamination $\mu$,  as constructed 
in
\cite{Thu}\cite{Ker}\cite{EpsMar}. 
We consider the
\emph{infinitesimal left
earthquake vector} $e^\mu_m=\frac
d{dt}E^{t\mu}_m{}_{|t=0} \in T_m
\mathcal T(S)$. This provides a
section
$e^\mu : \mathcal T(S)
\rightarrow T\mathcal T(S)$ of
the tangent bundle of $\mathcal
T(S)$, defined by $m\mapsto
e^\mu_m$. 

There similarly exists a
\emph{right earthquake}
$E^{-\mu}_m$ along $\mu$; the
notation is justified by the fact
that $m\mapsto E^{-\mu}_m$ is
the inverse of  $m\mapsto
E^{\mu}_m$. We can then consider
the \emph{infinitesimal right
earthquake vector}
$e^{-\mu}_m=\frac
d{dt}E^{-t\mu}_m{}_{|t=0}
\in T_m \mathcal T(S)$ and the
corresponding section
$e^{-\mu} : \mathcal T(S)
\rightarrow T\mathcal T(S)$ of
the tangent bundle of $\mathcal
T(S)$. Note that
$e^{-\mu}_m=-e^\mu_m$, but it
will be convenient to keep a
separate notation. 

Recall that the measured geodesic
laminations $\mu$, $\nu \in
\mathcal {ML}(S)$ \emph{fill up}
the surface $S$ if every
non-trivial measured geodesic
lamination has non-zero geometric
intersection number with at least
one of $\mu$, $\nu$. This is
equivalent to the condition that
every component of $S-\mu \cup
\nu$ is, either a topological disk
 bounded by the union of
finitely many geodesic arcs, or
a topological annulus bounded on
one side by the union of finitely
many geodesic arcs and going to a
cusp on the other side. 

\begin{prop} 
\label
{prop:TransEarthSectInTgtSpce}
Let $\mu$, $\nu \in \mathcal
{ML}(S)$ be two non-zero measured
geodesic laminations.  The
intersection of the two sections
$e^\mu$,
$e^{-\nu} : \mathcal T(S)
\rightarrow T\mathcal T(S)$ of
the tangent bundle of $\mathcal
T(S)$ is transverse. These
sections meet in exactly one
point if $\mu$ and $\nu$ fill up the
surface $S$, and are
otherwise disjoint. 
\end{prop}

\begin{proof}
We first translate the problem in
terms of the length functions
$l_\mu$, $l_\nu: \mathcal T(S)
\rightarrow
\mathbb R$ which to a metric
$m\in \mathcal T(S)$ associate
the $m$--lengths of the measured
geodesic lamination $\mu$ and
$\nu$. 
The Weil-Petersson
symplectic form on
$\mathcal T(S)$ induces an
isomorphism between its tangent
bundle $T\mathcal T(S)$ and its
cotangent bundle $T^*\mathcal
T(S)$. A celebrated result of
Scott Wolpert \cite{Wol} asserts
that this isomorphism sends the
section $e^\mu$ of $T\mathcal
T(S)$ to the section $dl_\mu$ of
$T^*\mathcal T(S)$. Therefore,
Proposition~\ref
{prop:TransEarthSectInTgtSpce} is
equivalent to showing that the
sections $dl_\mu$ and $-dl_\nu$
transversely meet in 1 or 0
point, according to whether $\mu$
and $\nu$ fill up the surface $S$ or
not. 

First consider the case where
$\mu$ and $\nu$ fill up $S$. The
intersection of the section
$dl_\mu$ and $-dl_\nu$ of
$T^*\mathcal T(S)$ correspond to
the points $m\in \mathcal T(S)$
where $d_ml_\mu=-d_ml_\nu$,
namely to the critical points of
the function $l_\mu +l_\nu
:\mathcal T(S) \rightarrow
\mathbb R$. It is proved in
\cite{Ker}\cite{Ker2} that,
because
$\mu$ and $\nu$ fill up $S$, the
function
$l_\mu+l_\nu$ admits a unique
critical point $m_0$. In
addition, the hessian of
$l_\mu+l_\nu$ at $m_0$ is positive
definite
\cite{Wol2}\cite{Ker2}. 

Let $u=d_{m_0}l_\mu =
-d_{m_0}l_\nu \in T^*_{m_0}
\mathcal T(S)$.
The intersection of the tangent
spaces of the sections
$dl_\mu$ and $-dl_\nu$
at $u$ consists of all the
vectors of the form
$T_{m_0}\left( dl_\mu\right)(v) = 
T_{m_0}\left( -dl_\nu \right)(v)$
for some $v \in T_{m_0}\mathcal
T(S)$, where $T_{m_0}\left(
dl_\mu\right)$, $T_{m_0}\left(
-dl_\nu \right):
T_{m_0}\mathcal T(S)
\rightarrow T_u T^* \mathcal
T(S)$ denote the tangent maps of
the sections $
dl_\mu$, $
-dl_\nu :
\mathcal T(S)
\rightarrow T^* \mathcal
T(S)$. 

Note that, for an
arbitrary $w \in T_{m_0}\mathcal
T(S)$, the vectors 
$T_{m_0}\left( dl_\mu\right)(w)$
and $ T_{m_0}\left( -dl_\nu
\right)(w)$ both project to $w$
through the tangent map of the
projection $T^*\mathcal T(S)
\rightarrow \mathcal T(S)$. In
particular, the difference 
$T_{m_0}\left( dl_\mu\right)(w)-
T_{m_0}\left( -dl_\nu
\right)(w)$ is tangent to the
fiber, namely is an element of
$T^*_{m_0}\mathcal T(S) \subset
T_u T^* \mathcal T(S)$. Using
local coordinates, one easily
sees that this difference is
 the image of $w$ under
the homomorphism $T_{m_0} \mathcal
T(S) \rightarrow T^*_{m_0}
\mathcal T(S)$ induced by the
hessian of the function $l_\mu +
l_\nu$. (Beware of optimistic
simplifications, though: for
instance, $T_{m_0}\left( -dl_\nu
\right)$ is not the same as
$-T_{m_0}\left(dl_\nu
\right)$ since they take their
values in different tangent
spaces.) Since the hessian is
non-degenerate, this homomorphism
is actually an isomorphism.
	
A consequence of this analysis is
that, if $T_{m_0}\left( dl_\mu\right)(v) = 
T_{m_0}\left( -dl_\nu
\right)(v)$, then $v$ must be 0.
Therefore, the intersection of
the two sections $dl_\mu$ and
$-dl_\nu$ of $T^* \mathcal T(S)$
at the point
$u$ is transverse, since these
sections have dimension $\theta$
and the total space $T^* \mathcal
T(S)$ has dimension $2\theta$. By
Weil-Petersson duality, this
proves that the sections
$e^\mu$ of $e^{-\nu}$ of
$T\mathcal T(S)$ have a
transverse intersection
consisting of exactly one point.

We now consider the case where
$\mu$ and $\nu $ do not fill up the
surface. In this case, the
function $l_\mu
+l_\nu$ has no critical point
\cite{Ker2} (see also the
discussion in \cite[\S
4.2]{Ser1}). By the above
translation, it follows that the
sections $e^\mu$,
$e^{-\nu}$ are disjoint.

This concludes the proof of
Proposition~\ref
{prop:TransEarthSectInTgtSpce}.
\end{proof}

If the measured geodesic
laminations $\mu$, $\nu$ fill up the
surface $S$, let $\kappa(\mu,
\nu)$ denote the (unique) critical
point of the length function
$l_\mu +l_\nu :\mathcal T(S)
\rightarrow
\mathbb R$. As above,
$\kappa(\mu, \nu )$ is also the
unique $m\in \mathcal T(S)$ such
that $e^\mu_m = e^{-\nu}_m$.

\begin{lem}
\label{lem:UniqueInfEarth}
If $\kappa(\mu', \nu)=
\kappa(\mu, \nu)$, then
$\mu'=\mu$. 
\end{lem}

\begin{proof} By an
(easy) infinitesimal version of
\cite{Thu} (see also \cite{Ker2}),
an infinitesimal left earthquake
completely determines the
measured geodesic lamination
along which it is performed. If
$e^{\mu'}_m = e^{-\nu}_m=
e^\mu_m$, it follows that
$\mu'=\mu$.
\end{proof}

We will need a result similar to
Proposition~\ref
{prop:TransEarthSectInTgtSpce}
in the unit
tangent bundle
$T^1\mathcal T(S)$. Recall that
the fibers $T_m^1\mathcal T(S)$ of
this bundle are the  quotient of
$T_m\mathcal T(S)-\{0\}$ under
the equivalence relation which
identifies $v$ to $tv$ when $t\in
\left]0,\infty\right[$. In
particular, $T^1\mathcal T(S)$ is
a manifold of dimension
$2\theta-1$.  Let
$\overline e{}^\mu$, $\overline
e{}^{-\nu}: \mathcal T(S)
\rightarrow T^1 \mathcal T(S)$ be
the sections induced by $e^\mu$
and $e^{-\nu}$. 

\begin{prop}
\label
{prop:TransEarthSectInUnitTgtSpce}
Let $\mu$, $\nu \in \mathcal
{ML}(S)$ be two non-zero measured
geodesic laminations.  The
intersection of the two sections
$\overline e{}^\mu$, $\overline
e{}^{-\nu}: \mathcal T(S)
\rightarrow T^1 \mathcal T(S)$  of
the unit tangent bundle of
$\mathcal T(S)$ is transverse.
If $\mu$ and $\nu$ fill up the
surface $S$, these sections meet
along a section above a line
$K(\mu, \nu)$ properly embedded
in $\mathcal T(S)$. If $\mu$,
$\nu$ do not fill up $S$, the
intersection is empty.
\end{prop}

\begin{proof}The two sections
meet above $m \in \mathcal T(S)$
if 
$\overline e{}^\mu_m=\overline
e{}^{-\nu}_m$, namely if there is
a
$t>0$ such that $e^{-\nu}_m=
te^\mu_m = e^{t\mu}_m$. By
Proposition~\ref
{prop:TransEarthSectInTgtSpce},
this can occur only when $t\mu$
and $\nu$ fill the surface,
namely only when $\mu$ and $\nu$
fill the surface. Consequently,
the two sections have empty
intersection if $\mu$ and $\nu$
do not fill up the surface.

If $\mu$ and $\nu$ fill up the
surface then, for every $t>0$,
Proposition~\ref
{prop:TransEarthSectInTgtSpce}
shows that there is a unique
$m=\kappa(t\mu, \nu) \in
\mathcal T(S)$ such that
$e^{t\mu}_m = e^{-\nu}_m$. As a
consequence, the two sections
 $\overline e{}^\mu$, $\overline
e{}^{-\nu}: \mathcal T(S)
\rightarrow T^1 \mathcal T(S)$
meet exactly above the image
$K(\mu, \nu)$ of the map
$\left]0,\infty\right[
\rightarrow \mathcal T(S)$
defined by $t\mapsto \kappa
(t\mu, \nu)$.

If $m=\kappa(t\mu, \nu)$ so that
$e^{t\mu}_m = e^{-\nu}_m$, the
tangent space
$T_{e^{t\mu}_m}T\mathcal T(S)$ is
the sum of the tangent spaces $T_m
e^{t\mu} \left ( T_m \mathcal T(S)
\right)$ and $T_m
e^{-\nu} \left ( T_m \mathcal T(S)
\right)$ of the sections
$e^{t\mu}$ and $e^{-\nu}$, by the
transversality property of
Proposition~\ref
{prop:TransEarthSectInTgtSpce}.
Therefore, in the unit tangent
bundle, the tangent space
$T_{\overline
e{}^{\mu}_m}T^1\mathcal T(S)$ is
the sum of $T_m \overline
e{}^{\mu}
\left ( T_m \mathcal T(S)
\right)$ and $T_m \overline
e{}^{-\nu} \left ( T_m \mathcal
T(S)
\right)$. As a consequence, the
intersection of the two sections
$\overline e{}^\mu$ and $\overline
e{}^{-\nu}$ is transverse above
the point  $m=\kappa(t\mu, \nu)$. 

By transversality, the
intersection of the sections is a
submanifold of the image of $\overline
e{}^\mu$. Its dimension is equal
to 1, by consideration of the
dimensions of $\mathcal
T(S)$ and $T^1\mathcal T(S)$.
Since the projection
$\overline e{}^\mu\left( \mathcal
T(S) \right) \rightarrow \mathcal
T(S)$ is a diffeomorphism, it
follows that the projection
$K(\mu, \nu)$ of this projection
is a 1--dimensional submanifold
of
$\mathcal T(S)$. 

By definition, $\kappa(t\mu,
\nu)$ is the unique minimum of the
convex function
$tdl_\mu+dl_\nu$, which has
positive hessian at this minimum.
It follows that $\kappa(t\mu,
\nu)$ is a continuous function of
$t$. Conversely, $t$ is
completely determined by
$m=\kappa(t\mu, \nu)$ by
Lemma~\ref
{lem:UniqueInfEarth}. As a
consequence, if
$m=\kappa(t\mu, \nu)$ stays in a
bounded subset of $\mathcal
T(S)$, then $t$ stays in a
compact subset of $\left]0,\infty
\right[$. In other words, the map
$\left]0,\infty\right[
\rightarrow \mathcal T(S)$
defined by $t\mapsto \kappa
(t\mu, \nu)$ is injective,
continuous and proper. It follows
that its image $K(\mu, \nu)$,
which we already know is a
1--dimensional submanifold of
$\mathcal T(S)$, is a line
properly embedded in
$\mathcal T(S)$. 
\end{proof}

Following the terminology of
\cite{Ser1} (motivated by
\cite{Ker2}), let the
\emph{Kerckhoff line} be the
proper 1--dimensional submanifold
$K(\mu,
\nu)
\subset
\mathcal T(S)$, consisting of all
the $\kappa(t\mu, \nu)$ with
$t>0$.

\section{Necessary condition for
small bending}
\label{sect:NecCond}

Let $t\mapsto m_t$, $t\in
\left[0, \varepsilon \right[$,  be
a small differentiable curve in
$\mathcal Q(S)$. If $\beta (m_t)
= \left( \beta^+(t),
\beta^-(t)\right) \in \mathcal
{ML}(S)^2$ is its bending measured
geodesic lamination, it is shown
in \cite{Bon98} that the right
derivative
$\frac d{dt^+} \beta(m_t)_{|t=0}$
exists, as an element of the
tangent space of $\mathcal
{ML}(S)^2$ at $\beta(m_0)$. In
general, because $\mathcal
{ML}(S)$ is not a differentiable
manifold, this tangent space
consists of geodesic laminations
with a transverse
structure which is less regular
than a transverse measure
\cite{Bon97}. However, if we
assume in addition that the
starting point
$m_0$ of the curve is fuchsian,
the tangent space of $\mathcal
{ML}(S)^2$ at $\beta(m_0)=(0,0)$
is just $\mathcal
{ML}(S)^2$; see \cite{Bon97}. 

We can therefore consider the
converse problem: Given a
fuchsian metric $m_0 \in
\mathcal F(S)$ and a pair $(\mu,
\nu) \in \mathcal {ML}(S)^2$ of
measured geodesic laminations,
does there exist a small
differentiable curve $t\mapsto
m_t\in \mathcal Q(S)$, $t\in
\left[0, \varepsilon \right[$,
originating from $m_0$ and such
that $\frac
d{dt^+} \beta(m_t)_{|t=0}=(\mu,
\nu)$? The following result shows
that $m_0$ is completely
determined by
$\mu$ and
$\nu$.

Note that, by construction, there
is a natural identification
between the submanifold $\mathcal
F(S)
\subset \mathcal Q(S)$
consisting of all fuchsian
metrics and the Teichm\"uller
space $\mathcal T(S)$.

\begin{prop}
\label{prop:NecCondition}
Let $\mu$, $\nu \in \mathcal
{ML}(S)$ be two measured
geodesic laminations, and let $t
\mapsto m_t$,
$t\in  \left[ 0,\varepsilon
\right[$ be a differentiable
curve in $\mathcal Q(S)$,
originating from a fuchsian
metric $m_0$ and such that the
derivative $\frac d{dt^+}
\beta(m_t)_{|t=0}$ of the bending
measured lamination is equal to
$(\mu, \nu)$. Then
$\mu$ and $\nu$ fill up the surface
$S$, and  $m_0
\in
\mathcal F(S) = \mathcal T(S)$ is
equal to the  minimum
$\kappa(\mu,
\nu)$ of the length function 
$l_\mu+l_\nu:\mathcal T(S)
\rightarrow \mathbb R$. 
\end{prop}

\begin{proof}We consider two
other curves in $\mathcal Q(S)$.

The first one is the \emph{pure
bending} $t\mapsto
B^{t\mu}_{m_0}$, obtained by
bending the surface $S$ along the
measured geodesic lamination
$t\mu$ while keeping the metric
induced on this pleated surface
equal to $m_0$. For
$t\geq 0$ small enough,
$B^{t\mu}_{m_0}$ is a
quasifuchsian metric for which
the positive side of the
boundary $\partial
C(B^{t\mu}_{m_0})$ is a pleated
surface with induced metric $m_0$
and with bending measured
geodesic lamination $t\mu$. See
\cite[\S 3]{EpsMar} or
\cite{Bon96} for the construction
of $B^{t\mu}_{m_0}$, and \cite[\S
9]{Mar} to guarantee that it is
quasifuchsian for $t$
sufficiently small. In addition,
it is proved in \cite[\S
3.9]{EpsMar}
\cite{Bon96} that this curve is
differentiable in
$\mathcal Q(S)$, and in
particular admits a tangent
vector $b^{\mu}_{m_0} = \frac
d{dt} B^{t\mu}_{m_0}{}_{|t=0} \in
T_{m_0} \mathcal Q(S)$. This
tangent vector $b^\mu_{m_0}$ is
the \emph{infinitesimal pure
bending} of $m_0\in \mathcal
F(S)$ along the measured geodesic
lamination $\mu$. 

The second curve will use 
the shear-bend
coordinates associated, as in 
\cite{Bon96}, to a maximal
geodesic lamination $\lambda$
containing the support of
$\mu$. These coordinates provide a
local parametrization of
$\mathcal Q(S)$ in terms of the
geometry of a pleated surface
with pleating locus $\lambda$.
Let $m'_t$ correspond to a
pleated surface whose induced
metric is equal to the
metric $m^+_t \in
\mathcal T(S)$ induced on the
positive component of the boundary
$\partial C(m_t)$ of the convex
core, and whose bending data is
equal to $t\mu$. By
\cite[\S9]{Mar},
$m_t$ is a quasifuchsian metric
for $t$ small.

It is proved in \cite{Bon98}
that, because the curve
$t\mapsto m_t$ is differentiable
and because $\frac d{dt^+}
\beta^+(m_t)_{|t=0}=\mu$, the
right derivative $\dot m^+_0 =
\frac d{dt} m_t^+{}_{|t=0} \in
T_{m_0} \mathcal T(S)$ exists and
the two curves
$t\mapsto m_t$ and $t\mapsto
m_t'$ have the same tangent
vector at $t=0$. In particular,
by differentiability of the
shear-bend coordinates,  the
tangent vector $\dot m_0 =
\frac d{dt} m_t{}_{|t=0}=
\frac d{dt} m_t'{}_{|t=0}\in
T_{m_0} \mathcal Q(S)$ is the sum
of $b^\mu_{m_0}$ and of 
$\dot m^+_0 =
\frac d{dt} m_t^+{}_{|t=0}\in
T_{m_0} \mathcal T(S)
=T_{m_0} \mathcal F(S)$.

Similarly, bending $S$ in the
negative direction, we can
define the infinitesimal pure
bending vector 
$b^{-\nu}_{m_0} = -
b^{\nu}_{m_0}\in T_{m_0}
\mathcal Q(S)$ of $m_0$ along
$-\nu$. We can also consider the
metric $m^-_t \in \mathcal T(S)$
induced on the negative side of
$\partial C(m_t)$. Then, as
above, the vector $\dot m_0$ is
the sum of $b^{-\nu}_{m_0}$ and
of $\dot m^-_0 =
\frac d{dt} m_t^-{}_{|t=0}\in
T_{m_0} \mathcal T(S)
=T_{m_0} \mathcal F(S)$.

Finally, we will use the complex
structure of
$\mathcal Q(S)$ coming from the
fact that the isometry group of
$\mathbb H^3$ is
$\mathrm{PSL}_2(\mathbb C)$.
Indeed, considering the
holonomy of hyperbolic metrics
embeds $\mathcal Q(S)$ into the
space $\mathcal R(S)$ of
(conjugacy classes) of
representations
$\pi_1(S)
\rightarrow
\mathrm{PSL}_2(\mathbb C)$. This
representation space $\mathcal
R(S)$ is a complex manifold near
the image of
$\mathcal Q(S)$, and admits this
image as a complex submanifold
(an open subset if $S$ is
compact). The shear-bend local
coordinates are well behaved with
respect to this complex
structure; see \cite{Bon96}. In
particular, multiplication by
$\mathrm i=\sqrt{-1}$ exchanges
shearing and bending. We will use
the following two consequences of
this. First of all, at a fuchsian
metric $m_0$, the tangent space
$T_{m_0} \mathcal Q(S)$ is the
direct sum of $T_{m_0}\mathcal
F(S)$ and of $\mathrm i T_{m_0}
\mathcal F(S)$. In addition, the
infinitesimal pure bending vector
$b^\mu_{m_0}$ belongs to $\mathrm i T_{m_0}
\mathcal F(S)$ and is equal to
$\mathrm i e^\mu_{m_0}$, where
$e^\mu_{m_0} \in T_{m_0} \mathcal
T(S) = T_{m_0}\mathcal F(S)$ is
the infinitesimal earthquake
vector along $\mu$. 

Applying the decomposition 
$T_{m_0} \mathcal
Q(S)=T_{m_0}\mathcal F(S)
\oplus \mathrm i T_{m_0}
\mathcal F(S)$ to the vector 
$\dot m_0 = \dot m^+_0 +
b^\mu_{m_0} = \dot m^-_0 +
b^{-\nu}_{m_0}$, we conclude that
$b^\mu_{m_0} = b^{-\nu}_{m_0}$.
Multiplying by $-\mathrm i$, it
follows that $e^\mu_{m_0} =
e^{-\nu}_{m_0}$. As in the proof
of Proposition~\ref
{prop:TransEarthSectInTgtSpce},
this is equivalent to the
property that
$m_0=\kappa(\mu,
\nu)$. 
\end{proof}

\section{Realizing small bending}

The goal of this section is to 
prove a converse to
Proposition~\ref
{prop:NecCondition}, by
constructing in Proposition~\ref
{Prop:ConstrSmallBending} a
small curve of quasifuchsian
metrics
$t \mapsto m_t\in\mathcal Q(S)$,
$t\in
\left[0, \varepsilon \right[$,
such that
$\beta(m_t)=(t\mu, t\nu)$ for
every
$t$. 

 For the
measured geodesic lamination
$\mu
\in
\mathcal{ML}(S)$, let
$\mathcal P^+(\mu)$ (resp. 
$\mathcal P^-(\mu)$) be the
space of quasifuchsian metrics
$m$ such that the positive (resp.
negative) component of the convex
core boundary
$\partial C(m)$ has bending
measured geodesic lamination 
$t\mu$, for some
$t\in\left[0,\infty\right[$.

Recall that $\theta =
-3\chi(S)+p$ denotes the
dimension of the Teichm\"uller
space $\mathcal T(S)$.

\begin{lem}
\label{lem:PleatSubmfd}
The space $\mathcal P^\pm (\mu)$
is a submanifold-with-boundary of
$\mathcal Q(S)$, with dimension
$\theta+1$ and with boundary
$\mathcal F(S)$.
\end{lem}

\begin{proof} 
We can use the coordinates
developed in \cite{Bon96}, and
associated to a maximal geodesic
lamination
$\lambda$ containing the support
of $\mu$. These provide an open
differentiable embedding
$\varphi:\mathcal Q(S)
\rightarrow \mathcal T(S) \times
\mathcal H_0(\lambda; \mathbb
R/ 2\pi \mathbb Z)$. The first
component of
$\varphi(m)$ is the hyperbolic
metric induced on the unique
$m$--pleated surface
$f_m$ with pleating locus
$\lambda$. The second
component is the bending
transverse cocycle of $f_m$,
which belongs to the topological
group
$\mathcal H_0(\lambda; \mathbb
R/ 2\pi \mathbb Z)
\cong 
 \left( \mathbb R/ 2\pi \mathbb
Z \right)^\theta \oplus \mathbb
Z/2$ of all $ \left( \mathbb R/ 2\pi \mathbb
Z \right)$--valued
transverse cocycles for
$\lambda$ that satisfy a certain
cusp condition. In general, the
bending of the pleated surface
$f_m$ is measured by a transverse
cocycle and not by a measured
geodesic lamination because $f_m$
is not necessarily locally
convex. 

For notational convenience, it
is useful to lift $\varphi$ to
an embedding $\psi:\mathcal Q(S)
\rightarrow \mathcal T(S) \times
\mathcal H_0(\lambda; \mathbb
R)$, such that $\psi$ sends
$\mathcal F(S) \cong \mathcal
T(S)$ to $\mathcal T(S) \times
\{0\}$ by the identity. Here, 
$\mathcal H_0(\lambda; \mathbb
R) \cong \mathbb R^\theta$ denotes
the space of
$\mathbb R$--valued transverse
cocycle satisfying the cusp
condition. Such a
$\psi$ exists and is unique
because $\mathcal Q(S)$ is simply
connected.

The vector space
$\mathcal H_0(\lambda; \mathbb R)$
contains the transverse measure
(also denoted by $\mu$) of the
measured geodesic lamination
$\mu$, and therefore also
contains the two rays $
\left[0,\infty
\right[\mu$ and
$\left]-\infty, 0
\right]\mu$, consisting of all
positive (resp.
negative) real multiples of
$\mu$.

We claim that
$\mathcal
P^+(\mu)$ locally corresponds,
under $\psi$,  to the intersection
of
$\psi\left(\mathcal
Q(S)\right)$ with $\mathcal T(S)
\times \left[0,\infty
\right[\mu$. Clearly,
$\psi$ sends an element of
$\mathcal P^+(\mu)$ to $\mathcal T(S)
\times \left[0,\infty
\right[\mu$. Conversely, it
is proved in
\cite{KamTan} that the map 
$\eta: \mathcal Q(S) \rightarrow
\mathcal T(S) \times
\mathcal {ML}(S)$
which, to a
quasifuchsian metric $m \in
\mathcal Q(S)$, associates the
induced metric $m^+\in\mathcal
T(S)$ and the
bending measured
lamination $\beta^+(m)$ of the
positive boundary component
of the convex core $C(m)$ is
a local homeomorphism. In
addition, interpreting the ray
$\left[0,\infty
\right[\mu$ as a subset of both
$\mathcal H_0(\lambda; \mathbb
R)$ and $\mathcal{ML}(S)$, the
local inverse for
$\eta$ constructed in
\cite{KamTan} coincides with
$\psi$ on $\mathcal T(S)
\times \left[0,\infty
\right[\mu$. It follows that, if
$m\in\mathcal Q(S)$ is
sufficiently close to $m_0\in
\mathcal P^+(\mu)$ and if
$\psi(m) \in \mathcal T(S)
\times \left[0,\infty
\right[\mu$, the bending
measured lamination
$\beta^+(\mu)$ is equal
to the second component $t\mu$ of
$\psi(m)$. Therefore, a metric
$m$ near
$m_0\in
\mathcal P^+(\mu)$ is in $
\mathcal P^+(\mu)$ if and only
if $\psi(m) \in \mathcal T(S)
\times \left[0,\infty
\right[\mu$. 

A similar property holds for
$\mathcal P^-(\mu)$ by symmetry.
Therefore, under the
diffeomorphism $\psi$,
$\mathcal P^+(\mu)$, $\mathcal
P^-(\mu)$ and $\mathcal F(S)$
locally correspond to the
intersection of
$\psi\left(\mathcal
Q(S)\right)$ with $\mathcal T(S)
\times \left[0,\infty
\right[\mu$, $\mathcal T(S)
\times \left]-\infty, 0
\right]\mu$ and $\mathcal T(S)
\times \{0\}$, respectively.
\end{proof}

Given two measured geodesic
laminations $\mu$,
$\nu\in\mathcal{ML}(S)$, we
want to consider the
intersection of
$\mathcal P^+(\mu)$
and $\mathcal P^-(\nu)$. Note
that this intersection is far from
being transverse, since these two
$(\theta+1)$--dimensional
submanifolds  both contain the
$\theta$--dimensional submanifold
$\mathcal F(S)$ as their boundary.
For this reason, we will consider
the manifold-with-boundary $\check
{\mathcal Q}(S)$ obtained by
blowing-up ${\mathcal Q}(S)$ along
the submanifold $\mathcal F(S)$.
Namely,
$\check{\mathcal Q}(S)$ is the
union of $\mathcal Q(S)-\mathcal
F(S)$ and of the unit normal
bundle
$N^1\mathcal F(S)$ with the
appropriate topology.
Recall that the normal bundle
$N\mathcal F(S) \rightarrow
\mathcal F(S)$ is intrinsically
defined as the bundle whose fiber
$N_m\mathcal F(S)$ at
$m\in \mathcal F(S)$ is the
quotient $T_m\mathcal
Q(S)/T_m\mathcal F(S)$, and that
the fiber $N^1_m \mathcal F(S)$
of the unit normal bundle
$N^1\mathcal F(S) \rightarrow
\mathcal F(S)$ is the quotient of
$N_m\mathcal F(S)-\{0\}$ under the
multiplicative action of
$\mathbb R^+$. Considering a
tubular neighborhood of
$\mathcal F(S)$,
$\check{\mathcal Q}(S)$ is easily
endowed with a natural structure
of differentiable manifold with
boundary
$N^1\mathcal F(S)$. 

Exploiting the complex structure
of $\mathcal Q(S)$, identify the
normal bundle $N\mathcal F(S)$ to
$\mathrm i T\mathcal F(S)$. The
inclusion map
$\mathcal P^+(\mu)-\mathcal F(S)
\rightarrow
\mathcal Q(S)-\mathcal F(S)$
uniquely extends to an embedding
$\mathcal P^+(\mu)\rightarrow
\check
{\mathcal Q}(S)$, which to $m\in
\mathcal F(S)= \partial \mathcal
P^+(\mu)$ associates the unit
normal vector
$\overline b{}_m^\mu \in
N^1\mathcal F(S)$ which is in the
direction of the infinitesimal
pure bending vector
$b_m^\mu$ of $m$ along
$+\mu$. We can similarly define an
embedding
$\mathcal P^-(\nu)\rightarrow
\check
{\mathcal Q}(S)$, by associating
to
$m\in
\mathcal F(S)$ the unit normal
vector $-\overline b{}_m^\nu \in
N^1\mathcal F(S)$ which is in the
direction of the infinitesimal
pure bending vector $b_m^{-\nu}
=-b_m^\nu$ of $m$ along
$-\nu$. Let
$\check{\mathcal P}^+(\mu)$ and
$\check{\mathcal
P}^-(\nu)\subset\check{\mathcal
Q}(S)$ be the respective images of
these embeddings.

\begin{lem}
The subspaces $\check{\mathcal
P}^+(\mu)$ and $\check{\mathcal
P}^-(\nu)$ are
submanifolds of
$\check{\mathcal Q}(S)$, with
boundary contained in $\partial
\check{\mathcal Q}(S)$.
\end{lem}

\begin{proof} The
embedding
$\psi:\mathcal Q(S)
\rightarrow \mathcal T(S) \times
\mathcal H_0(\lambda; \mathbb R)$
of the proof of
Lemma~\ref{lem:PleatSubmfd} sends
the vector $b^\mu_m$ to the
vector $(0,\mu)$ in the tangent
space of $\mathcal T(S) \times
\mathcal H_0(\lambda; \mathbb R)$.
The result then immediately
follows from the fact that
$\psi$ locally identifies
$\mathcal
P^+(\mu)$, $\mathcal
P^-(\mu)$ and $\mathcal T(S)$
 to  $\mathcal T(S)
\times  \left[0,\infty
\right[\mu$, $\mathcal T(S)
\times \left]-\infty, 0
\right]\mu$ and $\mathcal T(S)
\times \{0\}$, respectively. This
proves Lemma~\ref
{lem:PleatSubmfd}.
\end{proof}

\begin{prop}
\label
{prop:TransInterInBdryBlowUp}
The boundaries
$\partial \check{\mathcal
P}^+(\mu)$ and $\partial
\check{\mathcal P}^-(\nu)$ have a
non-empty intersection if and
only if $\mu$ and $\nu$ fill up the
surface $S$. If $\mu$ and $\nu$
fill up $S$, the
intersection of
$\partial\check{\mathcal
P}^+(\mu)$ and
$\partial\check{\mathcal
P}^-(\nu)$  in $\partial\check{
\mathcal Q}(S)$ is transverse,
and is equal to the image of the
Kerckhoff line
$K(\mu, \nu)$ under the section
$m\mapsto
\overline b{}_m^\mu 
=-\overline b{}_m^\nu\in
\partial\check{\mathcal Q}(S)$.
\end{prop}
\begin{proof} The
diffeomorphism $T^1\mathcal
T(S)=T^1\mathcal F(S) 
\rightarrow 
\mathrm i T^1
\mathcal F(S)= N^1\mathcal
F(S) =\partial \check {\mathcal
Q}(S)$ defined by multiplication
by
$\mathrm i$ sends $\overline
e{}^\mu_m$ to $\overline
b{}^\mu_m$ and $\overline
e{}^{-\nu}_m$ to $\overline
b{}^{-\nu}_m$. This translates
Proposition~\ref
{prop:TransInterInBdryBlowUp} to
a simple rephrasing of
Proposition~\ref
{prop:TransEarthSectInUnitTgtSpce}.
\end{proof}

An immediate consequence of
Proposition~\ref
{prop:TransInterInBdryBlowUp} is
that the intersection of the two
$(\theta+1)$--dimensional
submanifolds $\check{\mathcal
P}^+(\mu)$ and $
\check{\mathcal P}^-(\nu)$ in
$\check{\mathcal Q}(S)$ is
transverse near the boundary 
$\partial \check{\mathcal Q}(S)$.
In particular, the intersection 
 $\check{\mathcal
P}^+(\mu) \cap
\check{\mathcal P}^-(\nu)$ is a
2--dimensional submanifold of
$\check {\mathcal Q}(S)$ near
$\partial\check {\mathcal
Q}(S)$, with boundary 
$\partial \check{\mathcal
P}^+(\mu) \cap \partial
\check{\mathcal P}^-(\nu)$
contained in $\partial\check 
{\mathcal
Q}(S)$. 

By definition, a metric $m \in
{\mathcal P}^+(\mu) \cap
{\mathcal P}^-(\nu)$ has bending
measured lamination
$\beta(m)=(t\mu, u\nu)$ for some
$t$, $u\geq 0$. This gives a
differentiable map
$\pi: {\mathcal P}^+(\mu) \cap
{\mathcal P}^-(\nu)
\rightarrow \mathbb R^2$, defined
by
$m
\mapsto (t,u)$. 

Let $\check{\mathbb
R}^2$ be obtained by blowing up
$\mathbb R^2$ along $\{0\}$.
Because $\pi^{-1}(0)=
\mathcal F(S)$, the map
$\pi$ lifts to a differentiable
map
$\check\pi:
\check{\mathcal P}^+(\mu) \cap
\check{\mathcal P}^-(\nu)
\rightarrow \check{\mathbb R}^2$.

\begin{lem} 
\label{lem:PiCheckLocalDiffeo}
The map $\check\pi:
\check{\mathcal P}^+(\mu) \cap
\check{\mathcal P}^-(\nu)
\rightarrow \check{\mathbb R}^2$
is a local diffeomorphism near
$\partial
\check{\mathcal P}^+(\mu) \cap
\partial
\check{\mathcal P}^-(\nu)$.
\end{lem}

\begin{proof} We will
prove that, at any
point
$p_0
\in \partial
\check{\mathcal P}^+(\mu) \cap
\partial
\check{\mathcal P}^-(\nu)$, the
(linear) tangent map
$T_{p_0}\check\pi: T_{p_0}
\check{\mathcal P}^+(\mu)
\cap
T_{p_0} \check{\mathcal P}^-(\nu)
\rightarrow
T_{\check\pi(p_0)}\check{\mathbb
R}^2$ is injective.

 Let $v\in 
T_{p_0} \check{\mathcal P}^+(\mu)
\cap
T_{p_0} \check{\mathcal
P}^-(\nu)$ be such that $T_{p_0}
\check\pi(v)=0$. Considering the
map ${\mathcal P}^+(\mu)
\rightarrow
\left[0,\infty\right[$ which to
$m\in \mathcal P^+(\mu)$
associates $\beta^+(m)/\mu\in
\left[0,\infty\right[$ and the
induced map $\check{\mathcal P}^+(\mu)
\rightarrow
\left[0,\infty\right[$, we see
that $v$ must necessarily be in
the tangent space of the boundary
$\partial\check{\mathcal
P}^+(\mu)$. Symmetrically, it
must be in
$T_{p_0} \partial{\check{\mathcal
P}}^-(\nu)$. Therefore, $v$ is
tangent to the intersection 
$\partial
\check{\mathcal P}^+(\mu) \cap
\partial
\check{\mathcal P}^-(\nu)$.

Let us analyze the restriction of
$\check \pi$ to the boundary 
$\partial
\check{\mathcal P}^+(\mu) \cap
\partial
\check{\mathcal P}^-(\nu)$. By
Proposition~\ref
{prop:TransInterInBdryBlowUp}, 
$\partial
\check{\mathcal P}^+(\mu) \cap
\partial
\check{\mathcal P}^-(\nu)$
is equal to the image of the
Kerckhoff line
$K(\mu, \nu)$ under the section
$m\mapsto
\overline b{}_m^\mu 
=-\overline b{}_m^\nu\in
\partial\check{\mathcal Q}(S)$.
Recall that an element of the
Kerckhoff line $K(\mu, \nu)$ is
of the form $m=\kappa(t\mu, \nu)$
for some $t>0$, which is
equivalent to the property that
$b{}^{t\mu}_m= -
b{}^{\nu}_m$. We will also need a
coordinate chart for
$\check{\mathbb R}^2$ near
$\check\pi(p_0)$. Noting that the
image of $\pi$ is contained in the
quadrant
$\left[0,\infty\right[^2$, we can
use for this the chart $\varphi:
\left]0,\infty\right[ \times
\left [0,\infty \right[
\rightarrow 
\check{\mathbb R}^2$ defined on
the interior by
$(x,y) \mapsto (xy, y)$.

We claim that, if $m=\kappa(t\mu,
\nu)\in K(\mu, \nu)$, then
$\varphi^{-1}\circ \check\pi
\left (\,
\overline b{}_m^\mu \right)$ is
just equal to $(t,0)$. To see
this, choose, in the
2--dimensional manifold 
$\check{\mathcal P}^+(\mu) \cap
\check{\mathcal P}^-(\nu)$, a
small curve $s\mapsto \check m_s$,
$s\in \left[0,\varepsilon
\right[$, such that $\check m_0 =
\overline b{}_m^\mu \in
\partial \check{\mathcal P}^+(\mu)
\cap \partial
\check{\mathcal P}^-(\nu)$ and
such that
$\frac d{ds^+}\check m_s{}_{|t=0}$
is not tangent to the boundary.
The curve
$s\mapsto
\check m_s$ projects to a
differentiable curve
$t\mapsto m_t \in \mathcal Q(S)$
with $m_0=m$. By definition of
${\mathcal P}^+(\mu)$ and
${\mathcal P}^-(\nu)$, the
bending measured lamination
$\beta(m_s)$ is of the form
$(t(s)\mu, u(s)\nu)$ for two
differentiable functions $t(s)$,
and $u(s)$ with $t(0)=u(0)=0$.
Since $\frac d{ds^+}\check
m_s{}_{|t=0}$ points away from
the boundary, the curve
$s\mapsto m_s$ is not tangent to
$\mathcal F(S)$ at $s=0$, and it
follows that at least one of the
derivatives $t'(0)$, $u'(0)$ is
non-trivial. If we apply
Proposition~\ref
{prop:NecCondition}, we conclude
that
$m=\kappa\left (t'(0)\mu,
u'(0)\nu\right)$. Since 
$m=\kappa\left (t\mu,
\nu\right)$, it follows that
$t'(0)/u'(0)=t$ by Lemma~\ref
{lem:UniqueInfEarth}. In
particular, $t(s)/u(s)$ tends to
$t$ as $s$ tends to 0. Therefore,
$\varphi^{-1}\circ \check\pi
\left (
\overline b{}_m^\mu \right)$,
which is the limit of 
$\varphi^{-1}\circ \check\pi
\left (
\check m_s \right) = 
\varphi^{-1}\circ
\pi
\left (
 m_s \right)= \varphi^{-1}\circ
\bigl(t(s), u(s)\bigr) = \bigl(
t(s)/u(s), u(s)\bigr)$ as $s>0$
tends to 0, is equal to $(t,0)$.

This computation shows that the
restriction of $\check\pi$ to the
boundary $\partial
\check{\mathcal P}^+(\mu) \cap
\partial
\check{\mathcal P}^-(\nu)$ is a
diffeomorphism onto its image. In
particular, if $v\in 
T_{p_0} \check{\partial\mathcal
P}^+(\mu)
\cap
T_{p_0} \check{\partial\mathcal
P}^-(\nu)$ is such that $T_{p_0}
\check\pi(v)=0$, then necessarily
$v=0$. 

This concludes the proof that the
tangent map 
$T_{p_0}\check\pi:
T_{p_0} \check{\mathcal P}^+(\mu)
\cap
T_{p_0} \check{\mathcal P}^-(\nu)
\rightarrow
T_{\check\pi(p_0)}\check{\mathbb
R}^2$ is injective. Since $\check
\pi$ sends the boundary of the
2-dimensional manifold
$
\check{\mathcal P}^+(\mu) \cap
\check{\mathcal P}^-(\nu)$ to the
boundary of the 2-dimensional
manifold $\check{\mathbb
R}^2$, this proves that 
$\check
\pi:
\check{\mathcal P}^+(\mu) \cap
\check{\mathcal P}^-(\nu)
\rightarrow\check{\mathbb R}^2$
is a local diffeomorphism near
$p_0 \in \partial
\check{\mathcal P}^+(\mu) \cap
\partial
\check{\mathcal P}^-(\nu)$.
\end{proof}

This immediately gives the
following converse to
Proposition~\ref
{prop:NecCondition}. 

\begin{prop}
\label{Prop:ConstrSmallBending}
Let $\mu$, $\nu \in \mathcal
{ML}(S)$ be two measured geodesic
laminations which fill up the
surface $S$, and let $m_0$ be the
minimum $\kappa(\mu, \nu)$ of the
length function $l_\mu +l_\nu$.
Then there is a small
differentiable curve $t\mapsto
m_t\in \mathcal Q(S)$, $t\in
\left[0,\varepsilon\right[$,
beginning at $m_0$ and such that
the bending measured lamination
$\beta(m_t) $ is equal to $(t\mu,
t\nu)$ for every
$t$. 
\end{prop}

\begin{proof} Consider the curve
$t\mapsto (t,t)$, $t \in \left[
0, \varepsilon \right[$, in
$\mathbb R^2$. By Lemma~\ref
{lem:PiCheckLocalDiffeo}, for
$\varepsilon$  small enough,
there is a curve $t\mapsto \check
m_t \in
\check{\mathcal P}^+(\mu) \cap
\check{\mathcal P}^-(\nu)$ such
that $t \mapsto \check\pi \left(
\check m_t\right)$ coincides with
the lift of $t\mapsto (t,t)$ to 
$\check{\mathbb R}^2$. By
definition of the map $\pi$, this
just means that the projection
$m_t\in \mathcal Q(S)$ of 
$\check m_t\in \check{\mathcal
Q}(S)$ is such that
$\beta(m_t)=(t\mu, t\nu)$. 
\end{proof}

\section{Parametrizing
quasi-fuchsian groups by their
small bending}

Recall that $\mathcal {FML}(S)$
denotes the open subset of
$\mathcal {ML}(S)^2$ consisting
of those pairs $(\mu, \nu)$ such
that
$\mu$ and $\nu$ fill up the surface
$S$.  

Let $\check{\mathcal F}
\mathcal{ML}(S)$ be obtained by
blowing up $\mathcal {FML}(S)
\cup \{(0,0)\}$  along
$\{(0,0)\}$. Namely, 
$\check{\mathcal F}
\mathcal{ML}(S)$ is formally
obtained from
$\mathcal {FML}(S)$ by extending
each ray
$\left]0,\infty\right[(\mu, \nu)$
to a semi-open ray
$\left[0,\infty\right[(\mu,
\nu)$, with the obvious topology.
Note that the boundary
$\partial\check{\mathcal F}
\mathcal{ML}(S)$ is just the
quotient space of $\mathcal
{FML}(S)$ under the
multiplicative action of
$\mathbb R^+$. 

For every $(\mu, \nu)\in
\mathcal {FML}(S)$,
Proposition~\ref
{Prop:ConstrSmallBending}
provides a maximal ray
$R_{\mu\nu}=
\left[0,\varepsilon_{\mu\nu}
\right[(\mu, \nu)$ in
$\mathcal{FML}(S) \cup
\{(0,0)\}$  and a differentiable
map
$\Phi_{\mu\nu}: R_{\mu\nu}
\rightarrow  Q(S)$ such
that $\Phi_{\mu\nu}(\mu', \nu')$
has bending measured lamination
$(\mu', \nu')$ for every $(\mu',
\nu') \in R_{\mu\nu}$ and such
that $\check{\mathcal P}^+(\mu)$
and $\check{\mathcal P}^-(\nu)$
meet transversely along
$\beta(R_{\mu\nu})$. Here, the
statement that
$R_{\mu\nu}$ is maximal means that
$\varepsilon_{\mu\nu} \in
\left]0,\infty\right]$ is maximal
for this property. 

Note that $R_{\mu\nu}$ and
$\varphi_{\mu\nu}$ depend only on
the orbit of $(\mu, \nu)$ under
the action of $\mathbb R^+$,
namely on the corresponding point
of $\partial \check{\mathcal
F}\mathcal{ML}(S)$. Let
$\check R_{\mu\nu}$ be the lift
of $R_{\mu\nu}$ in
$\check{\mathcal
F}\mathcal{ML}(S)$, and lift
$\Phi_{\mu\nu}$ to $\check
\Phi_{\mu\nu}: \check R_{\mu\nu}
\rightarrow
\check
{\mathcal Q}(S)$. In particular,
$\check \Phi_{\mu\nu}$ sends the
initial point of $\check
R_{\mu\nu}$ to the bending vector
$\overline b{}^\mu_m= \overline
b{}^{-\nu}_m \in N^1 \mathcal F(S)
=
\partial \check{\mathcal Q}(S)$
with $m=\kappa(\mu, \nu)$. 

Let $\check U\subset
\check{\mathcal F}\mathcal{ML}(S)$
denote the union of all the
$\check R_{\mu\nu}$, and let
$\check \Phi : \check U
\rightarrow \mathcal {FML}(S)$
restrict to $\check \Phi_{\mu\nu}$
on each $\check R_{\mu\nu}$. Note
that the property
that $\Phi_{\mu\nu}(\mu', \nu')$
has bending measured lamination
$(\mu', \nu')$ implies that the
$\check R_{\mu\nu}$ are pairwise
disjoint, so that $\check \Phi$ is
well-defined. We want to show
that $\check \Phi$ is continuous.

\begin{lem}
\label{lem:ContinuityPleatSbmfd}
As the measured geodesic
lamination $\mu$ tends to $\mu_0$
for the topology of $\mathcal
{ML}(S)$, the submanifold
$\mathcal P^\pm(\mu)$ tends
to $\mathcal P^\pm(\mu_0)$ for
the topology of $\mathrm C^\infty$
convergence on compact subsets.
\end{lem}
\begin{proof} We will use the
tools developed in \cite{Bon96}. 

Let $\mu_n\in\mathcal {ML}(S)$,
$n\in \mathbb N$, be a sequence
converging to $\mu_0$. Let
$\lambda_n$ be a maximal geodesic
lamination containing the support
of $\mu_n$. Passing to a
subsequence if necessary, we can
assume that, for the Hausdorff
topology, the geodesic lamination 
$\lambda_n$ converges to a
geodesic lamination $\lambda_0$,
which is necessarily maximal and
contains the support of $\mu_0$.

The shear-bend coordinates
associated to $\lambda_n$  provide
an open biholomorphic embedding
$\Phi_n: \mathcal Q(S)
\rightarrow \mathcal H_0\left(
\lambda_n;
\mathbb C/2\pi \mathrm i \mathbb
Z \right)$. Here $\mathcal H_0\left(
\lambda_n;
\mathbb C/2\pi \mathrm i \mathbb
Z \right)$ is the topological
group of 
$\mathbb C/2\pi \mathrm i \mathbb
Z $--valued transverse cocycles
for the maximal geodesic
lamination
$\lambda_n$ which satisfy the
cusp condition, and is isomorphic
to
$\left( \mathbb C/2\pi 
\mathrm i \mathbb
Z \right)^\theta \oplus \mathbb
Z/2$. For a metric
$m\in
\mathcal Q(S)$, the real part of
$\Phi_n(m) \in \mathcal H_0\left(
\lambda_n;
\mathbb C/2\pi \mathrm i \mathbb
Z \right) =\mathcal H_0\left(
\lambda_n;
\mathbb R \right) \oplus \mathrm i
\mathcal H_0\left(
\lambda_n;
\mathbb R/2\pi  \mathbb
Z \right)$ measures the induced
metric of the unique $m$--pleated
surface with bending locus
$\lambda_n$, and the imaginary
part measures its bending. In
particular, $\mathcal P^\pm
(\mu_n)$ locally corresponds to
$\Phi_n^{-1} \left (
\mathcal H_0\left(
\lambda_n;
\mathbb R \right) \oplus \mathrm i
\mathbb R\mu_n
\right)$, or more precisely to a
branch of the immersion which is
the composition of the projection
$\mathcal H_0\left(
\lambda_n;
\mathbb R \right) \oplus \mathrm i
\mathbb R\mu_n \rightarrow
\mathcal H_0\left(
\lambda_n;
\mathbb C/2\pi \mathrm i \mathbb
Z \right)$ and of $\Phi_n^{-1}$.

To compare the various 
$\mathcal H_0\left(
\lambda_n;
\mathbb C/2\pi \mathrm i \mathbb
Z \right)$ pick a train track
$\tau$ carrying $\lambda_0$.
Since $\lambda_n$ converges to
$\lambda_0$ for the Hausdorff
topology, $\tau$ also carries the
$\lambda_n$ for $n$ large enough.
Then there is a well-defined
isomorphism $\Psi_n: \mathcal
H_0\left(
\lambda_n;
\mathbb C/2\pi \mathrm i \mathbb
Z \right) \rightarrow \mathcal
W\left(
\tau;
\mathbb C/2\pi \mathrm i \mathbb
Z \right)$, where $ \mathcal
W\left(
\tau;
\mathbb C/2\pi \mathrm i \mathbb
Z \right)$ is the group of $\mathbb C/2\pi \mathrm i \mathbb
Z $--valued edge
weights for $\tau$ satisfying the
switch and cusp relations. 

Because $\lambda_n$ converges to
$\lambda_0$ for the Hausdorff
topology, if follows from the
explicit construction of \cite[\S
5, \S 8]{Bon96} that
$\Phi_n^{-1}\circ
\Psi_n^{-1}$ converges to
$\Phi_0^{-1}\circ
\Psi_0^{-1}$, uniformly on
compact subsets of the image of
$\Psi_0 \circ \Phi_0$. Because
these maps are holomorphic, the
convergence is actually $\mathrm
C^\infty$.  Since $\mu_n$
converges to
$\mu_0$ for the topology of
$\mathcal {ML}(S)$, the edge
weight system $\Psi_n(\mu_n)$
converges to $\Psi_0(\mu_0)$ in
$\mathcal W_0(\tau, \mathbb R)$.
It
follows that 
$\mathcal P^\pm
(\mu_n)$, which locally
corresponds to
$\Phi_n^{-1} \circ
\Psi_n^{-1}\left (
\mathcal W_0\left(
\tau;
\mathbb R \right) \oplus \mathrm i
\mathbb R \Psi_n(\mu_n)
\right)$, converges to 
$\mathcal P^\pm
(\mu_0)$, which locally
corresponds to
$\Phi_0^{-1} \circ
\Psi_0^{-1}\left (
\mathcal W_0\left(
\tau;
\mathbb R \right) \oplus \mathrm i
\mathbb R \Psi_0(\mu_0)
\right)$, in the topology of $\mathrm
C^\infty$--convergence on compact
subsets.
\end{proof}

Note that we actually proved real
analytic convergence in
Lemma~\ref
{lem:ContinuityPleatSbmfd}.
However, we will only need
$\mathrm C^2$ convergence. 

\begin{thm} 
\label {them:HomeoInCheck}
The subset $\check U$
is an open neighborhood of
$\partial \check{\mathcal
F}\mathcal{ML}(S)$ in
$\check{\mathcal
F}\mathcal{ML}(S)$, and
$\check\Phi$ is a homeomorphism
from $\check U$ to an open
neighborhood $\check V$ of
$\partial \check{\mathcal Q}(S)$
in $\check{\mathcal Q}(S)$. 
\end{thm}

\begin{proof} The restriction
$\check \Phi_{\mu\nu}$ of $\check
\Phi$ to $\check R_{\mu\nu}$ was
constructed by considering the
transverse intersection of the
submanifolds $\check {\mathcal
P}^+(\mu)$ and $\check {\mathcal
P}^-(\nu)$ near the boundary of
$\check {\mathcal Q}(S)$. By
Lemma~\ref
{lem:ContinuityPleatSbmfd}, 
$\check {\mathcal
P}^+(\mu)$ and $\check {\mathcal
P}^-(\nu)$ depend continuously on
$(\mu, \nu)$ for the topology of
$\mathrm C^1$ convergence. (Note
that one needs the $\mathrm C^2$
continuity of $ {\mathcal
P}^+(\mu)$ and ${\mathcal
P}^-(\nu)$ to guarantee the
$\mathrm C^1$ continuity of
$\check {\mathcal P}^+(\mu)$ and
$\check {\mathcal P}^-(\nu)$ near
the boundary.) It follows that
the length
$\varepsilon_{\mu\nu}$ of
$R_{\mu\nu}=\left [0,
\varepsilon_{\mu\nu}
\right[(\mu, \nu)$ is a lower
semi-continuous function of
$(\mu, \nu)$, and that
$\Phi_{\mu\nu}$ depends
continuously on $(\mu, \nu)$.
This proves that the union $\check
U$ of the $\check R_{\mu\nu}$ is
open in
$\partial
\check{\mathcal
F}\mathcal{ML}(S)$, and that
$\check \Phi$ is continuous. 

Let $p:\check {\mathcal Q}(S)
\rightarrow \mathcal Q(S)$ be the
natural projection and, as usual,
let $\beta: \mathcal Q(S)
\rightarrow \mathcal{ML}(S)^2$ be
the bending map. By construction,
$\beta \circ p \circ \check \Phi$
is the identity map on $\check
U- \partial \check{\mathcal
F}\mathcal{ML}(S) \subset
\mathcal {ML}(S)^2$. It follows
that
$\check \Phi$ is injective on 
$\check U- \partial
\check{\mathcal
F}\mathcal{ML}(S)$, and therefore
on all of $\check U$ by
Proposition~\ref
{prop:NecCondition}. 

The two spaces
$\check U \subset \check{\mathcal
F}\mathcal{ML}(S)$ and
$\check{\mathcal Q}(S)$ are
topological
manifolds-with-boundary of the
same dimension $2\theta$. The map 
$\check \Phi : \check U
\rightarrow \mathcal {FML}(S)$ is
continuous and injective, and
sends boundary points to boundary
points. By the Theorem of
Invariance of the Domain, it
follows that its image $\check V
= \check \Phi\bigl ( \check U
\bigr)$ is open in $\check
{\mathcal Q}(S)$, and that
$\check \Phi$ restricts to a
homeomorphism $\check U
\rightarrow \check V$.
\end{proof}

\begin{thm} 
\label{thm:HomeoNoCheck}
There exists
an open neighborhood $V$ of the
fuchsian submanifold $\mathcal
F(S)$ in $\mathcal Q(S)$
such that the
bending map $\beta: \mathcal Q(S)
\rightarrow \mathcal{ML}(S)^2$
induces a homeomorphism between
$V-\mathcal F(S)$ and its image.  
\end{thm}
\begin{proof}
For the canonical identifications
between 
$\check {\mathcal Q}(S)-\partial
\check {\mathcal Q}(S)$ and $
{\mathcal Q}(S)- \mathcal F(S)$,
and between $\check{\mathcal
F}\mathcal{ML}(S)
-\partial \check{\mathcal
F}\mathcal{ML}(S)$ and
${\mathcal F}\mathcal{ML}(S)$,
the inverse of the restriction of
$\check\Phi$ to $\check U -
\partial
\check {\mathcal Q}(S)$ coincides
with $\beta$. Therefore, the only
part which requires some checking
is that the image $V$ of 
$\check V
= \check \Phi\bigl ( \check U
\bigr)$
 under the canonical projection 
$p:\check {\mathcal Q}(S)
\rightarrow \mathcal Q(S)$ is
open in $\mathcal Q(S)$. However,
this immediately follows from the
fact that the preimages of points
under $p$ are all compact, which
implies that $p$ is an open map.
(Note that this is false for the
projection $\check{\mathcal
F}\mathcal{ML}(S)
\rightarrow {\mathcal
F}\mathcal{ML}(S) \cup \{0\}$.)
\end{proof}

Theorem~\ref{thm:HomeoNoCheck} is
just Theorem~\ref{thm:MainThm}
stated in the introduction.
Theorem~\ref{thm:LocalImageBeta}
immediately follows from the
definition of
$\check U$ and from the fact that
$\check U$ is open
(Theorem~\ref{them:HomeoInCheck}).
\medskip 

We conclude with a few remarks. 

A recent result of Caroline Series
\cite{Ser2} shows that we can
restrict the neighborhood $V$ of
Theorem~\ref{thm:HomeoNoCheck}
so that
$\beta^{-1} \left(\beta(V)\right)
=V$.

 When $\mu$ and
$\nu$ are multicurves, namely
when their supports consist of
finitely many closed geodesics,
it follows from \cite{HodKer}
that the submanifolds 
$\check {\mathcal
P}^+(\mu)$ and $\check {\mathcal
P}^-(\nu)$ are everywhere
transverse, by a doubling
argument as in \cite{BonOta}.
Consequently, the open subset
$V$ of Theorem~\ref
{thm:HomeoNoCheck}
can be chosen so that the image
$U=\beta (V) \subset \mathcal
{ML}(S)^2$ contains all rays of
the form $\left[ 0, \infty \right[
(\mu, \nu)$ where $\mu$ and $\nu$
are multicurves. 

 As
indicated in the introduction, it
is conjectured that we can take
$V$ equal to the whole space
$\mathcal Q(S)$. See
\cite{BonOta} for a
characterization of the image of
$\mathcal Q(S)$ under $\beta$.


\begin{thebibliography}{EpM}


\bibitem[Bo1]{Bon96} 
Francis Bonahon,
\emph{Shearing hyperbolic 
surfaces, bending pleated
surfaces and {T}hurston's
symplectic form}, Ann. Fac. Sci.
Toulouse Math. (6)
\textbf{5}
(1996),  pp.~233--297.

\bibitem[Bo2]{Bon97}
\bysame,
\emph{Geodesic laminations with
transverse {H}\"older
distributions}, Ann. Sci. \'Ecole
Norm. Sup. (4)
\textbf{30} (1997),  pp.~205--240.
  

\bibitem[Bo3]{Bon98} 
\bysame, \emph{The boundary
geometry of $3$--dimensional
hyperbolic convex cores}, 
J. Differential Geometry
\textbf{50} (1998), pp.~1--23.


\bibitem[BoO]{BonOta}
Francis Bonahon, Jean-Pierre Otal,
\emph{Laminations mesur\'ees de 
plissage des vari\'et\'es
hyperboliques de dimension $3$},
preprint, 2001.

\bibitem[EpM]{EpsMar} 
 David B. A. Epstein, Albert
Marden, \emph{Convex hulls in
hyperbolic spaces, a theorem of
Sullivan, and measured pleated
surfaces}, in: \emph{Analytical
and geometric aspects of
hyperbolic space} (D.B.A.~Epstein
ed.),  L.M.S. Lecture Note Series
vol.
\textbf{111}, 1986, 
Cambridge University Press, 
pp.~113--253.

\bibitem[HoK]{HodKer} Craig D.
Hodgson, Steven P.
Kerckhoff, \emph {Rigidity of
hyperbolic cone-manifolds and
hyperbolic Dehn surgery}, J.
Differential Geom. \textbf{48}
(1998), pp.~1--59.

\bibitem[KaT]{KamTan} Yoshinobu
Kamishima, Ser P. Tan, 
\emph{Deformations spaces of
geometric structures}, in:
\emph{Aspects of low-dimensional
manifolds} (Y. Matsumo, S.
Morita eds.), Advanced
Studies  in Pure Math.
\textbf{20},
1992, Kinokuniya 
Company Ltd., Tokyo,
pp.  263--299.


\bibitem[Ke1]{Ker} 
Steven P. Kerckhoff,
 \emph{The Nielsen realization
problem},  Ann. of Math.
\textbf{117}  (1983), 
pp.~235--265. 

\bibitem[Ke2]{Ker2}
\bysame, 
\emph{Lines of minima in
Teichm\"uller space},   Duke Math.
J.  \textbf{65}  (1992),
pp.~187--213. 

\bibitem[Mar]{Mar} 
Albert Marden, \emph{The geometry 
of finitely generated Kleinian
groups},  Ann. of Math. \textbf 
{99} (1974), pp.~383--462.

\bibitem[Mas] {Mas} 
Bernard Maskit, 
\emph{Kleinian groups},
Grundlehren der Math.
Wiss. \textbf{287},
Springer-Verlag, 1988.

\bibitem[Se1]{Ser1} Caroline
M. Series,
\emph{On Kerckhoff minima and
pleating loci for quasi-Fuchsian
groups},  Geom. Dedicata 
\textbf{88}  (2001), 
pp.~211--237.

\bibitem[Se2]{Ser2} \bysame,
\emph{Limits of quasifuchsian
groups with small bending},
preprint, 2002.


\bibitem[Su]{Sul}
Dennis P. Sullivan,
\emph{Quasiconformal
homeomorphisms and dynamics II: 
Structural stability implies
hyperbolicity for  Kleinian
groups},  Acta Math. \textbf {155}
(1985), pp.~243--260.

\bibitem[Th1]{Thu1} 
William P. Thurston, \emph{The 
topology and geometry of
$3$--manifolds},  Lecture notes,
 Princeton University,  1976--79.

\bibitem[Th2]{Thu} 
\bysame, \emph{
Earthquakes in two-dimensional
hyperbolic geometry}, in: \emph{
Low-dimensional Topology and
Kleinian groups} 
(D.B.A.~Epstein ed.), 
L.M.S. Lecture Notes Series vol.
\textbf{112},  1986, 
Cambridge University press,
pp.~91--112.

\bibitem[Wo1]{Wol} 
Scott A. Wolpert, \emph{On
the symplectic geometry of
deformations of a hyperbolic
surface},  Ann. of Math.
\textbf{117}  (1983), 
pp.~207--234.

\bibitem[Wo2]{Wol2} 
\bysame, \emph{Geodesic 
length functions and the
Nielsen problem},  J. Differential
Geom.  \textbf{25} 
(1987), pp.~275--296.


\end{thebibliography}
\end{document}